\documentstyle[psfig,12pt]{amsart}
\newcommand{\C}{\mathbb C}

\newcommand{\Z}{\mathbb Z}
\renewcommand{\pf}{{\bf Proof. }}
\newcommand{\lam}{\lambda}

\newcommand{\sig}{\sigma}

\renewcommand{\epsilon}{\varepsilon}
\renewcommand{\phi}{\varphi}
\newtheorem{lem}{Lemma}[section]
\newtheorem{T}[lem]{Theorem}%[section]

\newtheorem{cor}[lem]{Corollary}%[section]
\newtheorem{defi}{Definition}[section]

\theoremstyle{definition}
\newtheorem{Rem}[lem]{Remark}%[section]
\newtheorem{Conjecture}{Conjecture}%[section]

\begin{document}

\title[ Braid Group Representations    ]{Irreducible Representations 
of Braid Groups of corank two      }
\author[I.Sysoeva        ]{Inna Sysoeva     }
\date{\today}
\address{Department of Mathematics \\
The Pennsylvania State University \\
University Park, PA 16802 }    
\email{sysoeva@@math.psu.edu   }

%%%%%%%%%%%%%%%%%%%

\begin{abstract}  This paper is the first part of a series of papers aimed at improving the classification by Formanek of the irreducible representations of Artin braid groups of small dimension. In this paper we classify all the irreducible complex representations $\rho$ of Artin braid group $B_n$ with the condition $rank (\rho (\sigma_i)-1)=2$  where $\sigma_i$ are the standard generators. For $n \geq 7$ they all belong to some one-parameter family of $n$-dimensional representations. 

\end{abstract}
\maketitle
%%%%%%%%%%%%%%%%%%%%%%%%%%%%%%%%%%%%%%%%%%%

%%%%%%%%%%%%%%%%%%%%%%%%%%%%%%%%%%%%%%%%%%%%%%%%%%%%%%%%
\section{Introduction.}

In his paper \cite{form} Edward Formanek classified all irreducible complex representations of Artin braid groups \hskip .3cm
$B_n$ \hskip .3cm
of dimension \hskip .3cm
at most \hskip .3cm
$n-1$.
This paper is the first in a series of papers aimed at extending this classification to irreducible representations of higher dimensions.

To describe our results, we need the following definition.

\begin{defi} The {\bf corank} of the representation
$\rho : B_n \to GL_r(\C)$ is $rank(\rho(\sig _i)-1)$ where the
$\sig_i$ are the standard generators of the group $B_n$ %{\rm (see \cite{chow}, %p.654).}

\end{defi}

\begin{Rem} Because the $\sig_i$ are conjugate to each other
(\cite{chow}, p.655), the number $rank(\rho(\sig _i)-1)$ does not depend on $i,$
which justifies the above definition.

\end{Rem}

The corank of  specializations of the reduced Burau representation (\cite{birm}, p.121; \cite{jones}, p.338)
and of the standard one-dimensional representation is $1.$ 

By the results of Formanek (\cite{form}, Theorem 23) almost all of the irreducible complex representations $B_n$ of degree at most $n-1$ of are the tensor product of
a one-dimensional representation and a representation of corank $1.$ He also classified all the irreducible representations of corank 1 (see \cite{form}, Theorem 10). For $n$ large enough they are one of the following.

\begin{enumerate}

\item A one-dimensional representation $\chi(y):B_{n} \to {\C}^*,\,\,\,\chi(y)(\sig_i)=y$

\item An irreducible $(n-1)-$dimensional specialization of the reduced Burau representation

\item An irreducible $(n-2)-$dimensional specialization of the composition factor
 of the reduced Burau representation

\end{enumerate}

The main goal of this paper is to classify all the irreducible complex representations of corank
$2.$ Apart from a number of exceptions for $n \leq 6,$ they
all are equivalent to  specializations for $u \neq 1,$ $u \in {\C}^*$
of the following representation $ \rho: B_n \to GL_n({\C }[u ^{\pm 1}]),$  first discovered by Dian-Ming Tong, Shan-De Yang and Zhong-Qi Ma in \cite{tym}:
$$\rho (\sig_i)=\left( \begin{array}{ccccc}
I_{i-1}&&&\\
&0&u&\\
&1&0&\\
&&&I_{n-1-i}
\end{array} \right),$$ 
\vskip .5cm
\noindent
for $i=1,2,\dots, n-1,$ where $I_{k}$ is the $k\times k$ identity matrix.

The main tool we use is the friendship graph of a representation. Namely the (full) friendship graph of  a representation $\rho$ of a braid group $B_n$ is a graph whose vertices are the set of generators $(\sigma_0,)$ $ \sigma_1, \dots , \sigma_{n-1}$ of $B_n$. Two vertices $\sigma_i$ and $\sigma_j$ are joined by an edge if and only if $Im(\rho (\sigma_i)-1) \cap Im(\rho (\sigma_j)-1) \neq \{0\}.$

Using the braid relations, we investigate the structure of the friendship graph. It turns out that  every irreducible representation of $B_n$ of dimension at least $n$ and corank $2$ the friendship graph is a chain, provided that $n\geq 6.$ This means that $\sigma_i$ and $\sigma_j$ are joined by an edge if and only if $|i-j|=1.$ 

For a given friendship graph it is relatively easy to classify all irreducible complex representations of $B_n$ for which it
is the associated friendship graph."
When the graph is a chain, we get specializations of the representation discovered by Tong, Yang and Ma.

Now we are going to explain the place of this paper in the coming series. According to \cite{form}, Theorem 23, for $n$ large enough every irreducible complex representation of $B_n$ of dimension at most $n-1$ is a tensor product of a one-dimensional representation and a representation of corank 1. Using similar ideas one can show that for $n$  large enough every  irreducible complex representation of $B_n$ of dimension at most $n$ is a tensor product of a one-dimensional representation and a representation of corank 2. Therefore one can use the results of this paper to extend the classification theorem of Formanek to the representations of $B_n$ of dimension $n.$ The proof of this result will appear elsewhere.

Another result, which will appear elsewhere is that for $n$ large enough there are no irreducible complex representations of $B_n$ of corank 3 and no irreducible complex representations of $B_n$ of dimension $n+1$. 

Based on the above result we would like to make the following two conjectures.

\begin{Conjecture}
For every $k\geq 3$ for $n$ large enough there are no irreducible complex representations of $B_n$ of corank $k.$
\end{Conjecture}

\begin{Conjecture}
For every $k\ge 1$ for $n$ large enough there are no irreducible complex representations of $B_n$ of dimension $n+k.$
\end{Conjecture}

We should also note that for the purpose of brevity we did not include in this paper some of the details of the classification of representations of $B_n$ for small $n.$ The full proof can be found in our thesis \cite{thesis}, Chapters 6 and 7.

The paper is organized as follows. In section 2 we introduce some convenient notation that will be used throughout the rest of the paper. In section 3
we define the friendship graph of the representation and study its structure. We also study the case when the friendship graph is totally disconnected. In section 4 we prove that for $n\ge 6$ for any irreducible complex representation of $B_n$ of corank $2$ and dimension at least $n$ the associated friendship graph is a chain. In section 5 we determine all irreducible representations of corank 2 whose friendship graph is a chain.

{\bf Acknowledgments:} The author would like to express her deep gratitude to professor Formanek for the numerous helpful discussions and comments on the preliminary versions of this paper, and for  generous financial support of this research.

%%%%%%%%%%%%%%%%%%%%%%%%%%%%%%%%%%%%%%%%%%%%%%%%%%
\section{Notation and preliminary results}

Let $B_n$ be the braid group on $n$ strings. It has a presentation
$$B_n=<\sig_1, \dots , \sig_{n-1}| \sig_i \sig_{i+1} \sig_i=\sig_{i+1} \sig_{i} \sig_{i+1}, 1 \leq i \leq n-2 ; \sig_i \sig_j = \sig_j \sig_i, |i-j| \geq 2>.$$

\begin{lem} \label{lem:21} For the braid group 
$B_n$
set $$\tau = \sig_1 \sig_2 \dots \sig_{n-1} {\textrm{ and }} \sig_0=\tau \sig_{n-1}\tau^{-1}.$$

Then:

1) {\rm (\cite{chow}, p.655)}  $$\sig_{i+1}=\tau \sig_{i} \tau^{-1},$$ for $1\leq i \leq n-2;$

2) $$\sig_i \sig_{i+1} \sig_i=\sig_{i+1} \sig_{i} \sig_{i+1},$$
$$\sig_{i+1}=\tau \sig_{i} \tau^{-1},$$
and $$ \sig_i \sig_j = \sig_j \sig_i, |i-j| \geq 2$$ for all $i,j$ where
indices are taken modulo $n.$
\end{lem}
\begin{Rem} Taking into account the above lemma, we  also have  the following presentation of $B_n:$
$$B_n=<\sig_0,\sig_1, \dots , \sig_{n-1}| \sig_i \sig_{i+1} \sig_i=\sig_{i+1} \sig_{i} \sig_{i+1};  \sig_i \sig_j = \sig_j \sig_i, |i-j| \geq 2; \sig_0=\tau \sig_{n-1} \tau ^{-1}>$$
for all $i,j$ where indices are taken modulo $n$ and $\tau$ is defined as above.

\end{Rem}
\vskip .5cm

Let $\rho : B_n \to GL_r (\C ) $ be a matrix representation of $B_n$
with $$\rho (\sig_i)=1+A_i,$$ and $$\rho (\tau)=T \in  GL_r (\C ) .$$
Then for any $i$ (indices are modulo $n$), the relation
$$\tau \sig_{i} \tau^{-1}=\sig_{i+1}$$ implies that $$TA_{i}T^{-1}=A_{i+1}.$$ 

Hence all the $A_i$ are conjugate to each other, so they have
the same rank, spectrum and Jordan normal form.

\begin{lem}\label{lem:2} For a representation $\rho$ of $B_n$ with  
$$\rho (\sig_i)=1+A_i,$$ we have:

1) $A_iA_j=A_jA_i,$ for $|i-j| \geq 2;$

2)$ A_i+A_i^2+A_iA_{i+1}A_i=A_{i+1}+A_{i+1}^2+A_{i+1}A_iA_{i+1}$ 

for
all  $i =0,1, \dots , n-1,$ where indices are taken modulo $n.$

\end{lem}
\pf 
 This  follows easily from the relations on the generators of $B_n.$

%%%%%%%%%%%%%%%%%%%%%%%%%%%%%%%%%%%%%%%%%%%%%%%%%%%%%%%%%%%%
\section {The friendship graph.}

In this section we define and prove some properties
of the {\it friendship graph} which is a finite graph
associated with a representation of $B_n.$ Our graphs are
{\it simple-edged,} which means that there is at most one
unoriented edge joining two vertices, and no edges joining 
a vertex to itself.

We assume throughout this section that we have a representation
$$\rho:B_n \to GL_r({\C}),$$ with
$$\rho(\sig_i)=1+A_i, \,\,\,\, (i = 0, 1, \dots ,n-1).$$

\begin{defi}

1) $A_i,\,\, A_{i+1}$ are 
{\bf neighbors} (indices modulo $n$).

2) $A_i,\,\, A_j$ are {\bf friends}
if $$ Im(A_i)\cap Im(A_j) \neq \{ 0 \}.$$

3) $A_i,\,\, A_j$ are {\bf  true friends}
if either

(a) $A_i$ and $A_j$  are not neighbors, and $$A_i  A_j =A_j  A_i\neq  0 ;$$  or 

(b) $A_i$ and $A_j$ are neighbors, and 

$$A_i +A_i^2+A_iA_j  A_i=A_j +A_j^2+A_j A_i  A_j \neq 0 .$$

\end{defi}

\begin{lem} \label{lem:trfr} If $A,B$ are true friends, then they are friends.
\end{lem} 
\pf 1) If $A$ and $B$ are not neighbors, then $AB=BA\neq 0,$
so, $$Im(A) \cap Im(B) \supseteq Im(AB)\cap Im(BA)= Im(AB)  \neq \{0\}.$$

2) If $A$ and $B$ are  neighbors, then 
$$A(1+A+BA)=A+A^2+ABA=B+B^2+BAB=B(1+B+AB) \neq 0 ,$$ and again
$$Im(A) \cap Im(B) \supseteq Im(A+A^2+ABA) \neq \{0\}.$$ 
\vskip .5cm

%%%%%%%%%%%%%%%%%%%%%%%%%%%%%%

\begin{defi} {\bf The full friendship graph} (associated with
the representation 
$\rho:B_n \to GL_n({\C})$ )  is the simple-edged graph with
$n$ vertices $A_0, A_1, \dots , A_{n-1}$ and an edge joining
$A_i$ and $A_j$ ($i \neq j$) if and only if $A_i$ and $A_j$ are
friends.

 {\bf The friendship graph} is the subgraph with vertices
$A_1, \dots , A_{n-1}$ obtained from the full friendship graph
by deleting $A_0$ and all edges incident to it.
\end{defi}

%\vskip 1cm

%\psfig{figure=graph2.eps}
%\vskip .5cm
%Two vertices $A_i,A_j,$ $ i \neq j,$ are connected if and only if
 %they are friends:

%\vskip 1cm

%\psfig{figure=graph1.eps}

%\vskip .5cm

%\end{defi}

Our main interest is the friendship graph, but it is convenient
to introduce the full friendship graph as a tool, because of the following lemma.

\begin{lem} There is an edge between $A_i$ and $A_j$ in the
full friendship graph if and only if there is an edge between
$A_{i+k}$ and $A_{j+k}$ where indices are taken modulo $n$.
In other words, ${\Z }_n$ acts on the full friendship
graph by permuting the vertices cyclically.

\end{lem}
\pf This follows immediately  from the  fact that conjugation by
$T= \rho(\tau)=\rho(\sig_1 \dots \sig_{n-1})$ permutes $\sig_0, \sig_1, \dots , \sig_{n-1}$ cyclically (Lemma \ref{lem:21}).

\begin{lem}[Lemma about friends] \label{lem:fr} Let $A$ and $B$ be neighbors   which are not friends. If $C$ is not a neighbor
of $A$ and $C$ is a friend of $B$ then $C$ is a true friend of $A.$
\vskip 1.5cm

\psfig{figure=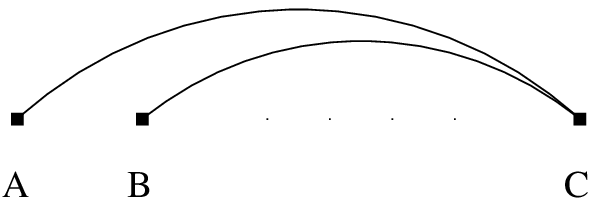}

\end{lem}
\vskip 1.5cm
\pf  By  lemma \ref {lem:trfr}, $A$ and $B$ are true  not friends, because   
they are not friends, that is
$$A+A^2+ABA=B+B^2+BAB=0.$$

 Consider $y \in V$ such that $Cy \in Im(B), Cy=Bz \neq 0$ ($y$ exists because
$C$ and $B$ are friends).
Then $$BACy=BABz=-(B+B^2)z=-(1+B)Bz \neq 0$$ because $Bz\neq 0$ and
$(1+B)$ is invertible. 

So, $AC=CA \neq	0;$ that is, $A$ and $C$ are true friends.

\vskip .5cm

%%%%%%%%%%%%%%%%%%%%%%%%%%%%%%%%%%%%%%%%

\begin{T} \label {T:connect} Let $\rho:B_n \to GL_r({\C})$  be a representation.    Then one of the following holds.

(a) The full friendship graph is  totally disconnected
(no friends at all).

(b) The full friendship graph has an edge between $A_i$ and $A_{i+1}$ for all $i.$ 

(c)  The full friendship graph has an edge between $A_i$ and $A_{j}$ whenever
$A_i$ and $A_{j}$ are not neighbors.

\end{T}
\pf Suppose neither (a) nor (b) holds. Since the graph is not totally
disconnected, there is an edge joining some vertices
$B$ and $C.$ Since (b) does not hold, no neighbors are joined by an edge.
Lemma \ref{lem:fr} implies that there is an edge between $C$ and any neighbor
of $B$ which is not a neighbor of $C.$ It follows inductively that there
is an edge joining $C$ to every vertex which is not a neighbor
of $C.$ Then (c) holds, because the full friendship graph is a ${\Z }_n$-graph. 

\vskip .5cm

\begin{defi} The friendship graph (the full friendship graph) is {\bf a chain},
if the only edges are between neighbors.

\end{defi}

Case $(b)$ of the above theorem can be restated as

{\it (b) The full friendship graph contains the chain graph.}

%For $n=2,n=3$ the statement of the theorem is
%true because any graph on 1 or 2 vertices is either connected
%or totally disconnected.

%Suppose now that $n \geq 5.$ Consider the full friendship graph of the %representation. Suppose that the graph is not totally disconnected -- that is, %there is an edge between 2 vertices
%$A_i$ and $A_j$ for some $i,j.$  Then we will prove the following statement:

 %There exist two neighbors $A_k$ and $A_{k+1}$ who are in the same connected %component.

%From this statement it follows easily that all the vertices of the full 
%friendship graph are
%in the same connected component by using the fact that the full friendship
%graph is a  ${\Z }_n$-graph. But if a ${\Z }_n$-graph is connected, then there 
%exists a covering path, starting from the vertex $A_0.$ (Covering path is a %path which passes through each vertex of a graph
%exactly once.) So, the
%friendship graph, obtained from the full friendship graph by deleting
 %the vertex $A_0$ and all edges
%adjacent to it, is also connected.

%Pick a vertex $A_i.$ Because  the full friendship graph is a ${\Z }_n$-graph,
%there is an edge between
%$A_i$ and $A_j$ for some $j.$ If $j=i+1$ then
%$A_i$ and $A_{i+1}$ are in the same connected component. Suppose now that
%$j \neq i+1.$ Then because $n \geq 5,$ at least one of $A_{i-1}$ and 
%$A_{i+1}$ is not a neighbor of $A_j.$ So, by the lemma about friends
%(\ref {lem:fr}), $A_{i-1}$ or $A_{i+1}$ ( or both) is in the same connected 
%component with $A_i.$

\vskip .5cm

\begin{cor} \label{cor:35} For $n \neq 4,$ the friendship graph and the full friendship graph
are either totally disconnected (no edges) or connected.

\end{cor}

%%%%%%%%%%%%%%%%%%%%%%%%%%%%%%%%%%%%%%%%%  

\begin{Rem} For $n=4$ there is a friendship graph which is neither totally disconnected nor connected:
\vskip 1cm
\psfig{figure=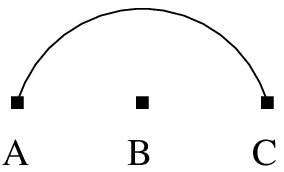}
\vskip 1cm
\end{Rem}

By \cite{thesis}, Lemmas 6.2 and 6.3, every  representation
of $B_4$ of corank $2$ and dimension at least $4,$ which has this friendship graph, is reducible.
\vskip .5cm
%%%%%%%%%%%%%%%%%%%%%%%%%%%%%%%%%%%%%%%

Now consider the case when the friendship graph is totally disconnected
(that is,  statement $(a)$ of theorem \ref{T:connect} holds).

\begin{lem} \label{lem:BB} If $A$ and $B$ are neighbors and not friends then:

(a) $A^2B=AB^2;$  $BA^2=B^2A.$

(b) If $x \in Im(A) \cap Ker(A- \lam I),$ then $B(Bx)=\lam (Bx)$ and
$ABx=-(1+\lam) x.$

\end{lem}
\pf (a). By  lemma \ref{lem:trfr}, $A$ and $B$ are  not true friends, so
$$A+A^2+ABA=B+B^2+BAB=0.$$
Multiplying the left hand side on the right by $B$ and the right hand side 
on the left by $A$ gives
$$AB+A^2B+ABAB=0=AB+AB^2+ABAB.$$
Thus, $A^2B=AB^2;$ by a symmetric argument $BA^2=B^2A.$

(b) Let $x=Ay \in  Im(A) \cap Ker(A- \lam I).$  
Then $$B(Bx)=B^2Ay=BA^2y=BAx=\lam Bx,$$ and
$$0=(A+A^2+ABA)y=(1+A+AB)x=(1+\lam) x+ABx.$$
Thus, $ABx=-(1+\lam) x.$

\vskip .5cm

%%%%%%%%%%%%%%%%%%%%%%%%%%%%%%%%%%

\begin{T} \label{T:38} Let  $\rho:B_n \to GL_r(\C),$ ($n \geq 2$) be an  irreducible representation,
 whose associated friendship graph
is totally disconnected. Then $r=dim V \leq n-1.$

\end{T}
\pf If $A_i = 0, $  $\rho$ is a trivial representation and $r =1.$

If $A_i\neq 0, $ 
choose an eigenvalue $\lam$ for $A_1$ and a non-zero vector $$x_1\in Im(A_1) \cap Ker(A_1-\lam I).$$  Set $x_2=A_2x_1, x_3=A_3x_2, \dots ,
x_{n-1}=A_{n-1}x_{n-2}, U=span \{x_1, x_2, \dots , x_{n-1} \}.$
By induction and lemma \ref{lem:BB} (b) $x_i\in Im(A_i) \cap Ker(A_i-\lam I).$

Let $x_i=A_iy_i.$ Then by  lemma \ref{lem:BB} (b) and the fact that
$A_iA_j=A_jA_i=0,$ if $i$ and $j$ are not neighbors,
$$A_{i-1}x_i=A_{i-1}A_{i}x_{i-1}=-(1+\lam)x_{i-1},\,\,\,\, i=2, \dots ,n-1,$$
$$A_{i}x_i=\lam x_i,\,\,\,\, i=1, \dots ,n-1,$$
$$A_{i+1}x_i=x_{i+1},\,\,\,\, i=1, \dots ,n-2,$$ and 
$$A_{j}x_i=A_{j}A_{i}y_{i}=0 \,\,\,\, j \neq i-1,i,i+1.$$

Thus $U$ is invariant under $B_{n}.$ Hence $r=dim U \leq n-1,$
since $\rho$ is irreducible. 

%%%%%%%%%%%%%%%%%%%%%%%%%%%%%%%%%%%%%%%%%%%%%%%%%%%%%%%%%%

\begin{cor} \label{cor:connect} Let $\rho:B_n \to GL_r(\C)$  be  irreducible, 
where $r=dim V \geq n,$ $n \neq 4.$ 

Then the associated friendship graph is connected.

\end{cor}
\pf By corollary \ref{cor:35} the friendship graph of $\rho$
is either totally disconnected or connected. By theorem \ref{T:38} it is not disconnected.

\begin{cor} \label{cor:dim} Let $\rho:B_n \to GL_r(\C)$  be  irreducible, 
where $r=dim V \geq n,$ $n \neq 4.$ Suppose $\rho (\sig_i)=1+A_i,$ where $rank (A_i)=k.$

Then $r=dim V \leq (n-1)(k-1)+1.$

In particular, for $k=2,$ $r=dim V = n,$ where $V={\C}^n.$

\end{cor}
%%%%%%%%%%%%%%%%%%%%%%%%%%%%%%%%%%%%%%%%%
\pf By  corollary \ref {cor:connect}, the friendship graph  of the representation  is connected. Arrange the vertices
of the graph in a sequence $A_{i_1}, A_{i_2},\dots , A_{i_{n-1}}$ such that each term $A_{i_j}, \,\,\, 2 \leq j \leq n-1,$
 is 
a friend of one the  terms $A_{i_1}, A_{i_2},\dots , A_{i_{j-1}}.$ Then

$$\dim(Im(A_{i_1}))=k$$
\vskip .2cm

$$\dim(Im(A_{i_1})+Im(A_{i_2}))\leq k+k-1=2k-1$$
$$ \dots$$
$$\dim(Im(A_{i_1})+ \dots +Im(A_{i_{n-1}}))\leq k+(n-2)(k-1)=(n-1)(k-1)+1.$$

\vskip .5cm
%%%%%%%%%%%%%%%%%%%%%%%%%%%%%%%%%%%%%%%%%%%%%%%%%%%%%%%%%%%%
Combining  Theorem \ref{T:connect} and Corollaries \ref{cor:connect} and
\ref{cor:dim}, we get the following

\begin{T} \label{T:dva} Let $\rho:B_n \to GL_r(\C)$  be  irreducible, 
where $r=dim V \geq n,$ $n \neq 4.$ Suppose $\rho (\sig_i)=1+A_i,$ where $rank (A_i)=2.$

Then $r =n$ and one of the following  holds.

(a) The full friendship graph has an edge between $A_i$ and $A_{i+1}$ for all $i.$ 

(b)  The full friendship graph has an edge between $A_i$ and $A_{j}$ whenever
$A_i$ and $A_{j}$ are not neighbors.

\end{T}

%From the above it follows that when $dim V \geq n,$ the friendship
%graph is connected. 
%Theorem \ref{T:discon} below shows that for $n\geq 6,$
%if $dim V \geq n$  then the friendship graph of an irreducible 
%representation contains a chain.
%
%\vskip .5cm

%%%%%%%%%%%%%%%%%%%%%%%%%%%%%%%%%%%%%%%%%%%%%%%%%%%%%%%%%%%%

\section{For corank 2 the friendship graph is a chain.}
%%%%%%%%%%%%%%%%%%%%%%%%%%%%%%%%%%%%%%%%%%%%%%%%%%%%%%%%%%%%%

In this section, we assume throughout that we have an irreducible representation 
$$\rho:B_n \to GL_r({\C}),$$ where $r\geq n,$ and 
$$\rho(\sig_i)=1+A_i,\,\,\,rank (A_i)=2,\,\,\,1 \leq i \leq n-1.$$

%\begin{lem} \label{lem:4gon} Let $\rho:B_n \to GL_r(\C)$  be a  
%representation, 
%where $n\geq 5,$ $r\geq n$ and $rank (A_1)=2.$ Suppose  that
 %there exist  distinct vertices $A,B,C,D$ such that $A,B$ are 
%friends, $B,C$ are friends,
%$C,D$ are friends and $D,A$ are friends.

 %Then $B,D$ are friends and $A,C$ are friends. In other words, if 
%there is a 4-gon of friends then it has both diagonals:

%\vskip 1cm
%\psfig{figure=pic4gon.eps}
%\vskip 1cm
%\end{lem}
%\pf  Suppose that $B$ and $D$ are not friends. Let $c_1 \in Im (C) 
%\cap Im(B)$
%and $c_2 \in Im (C) \cap Im(D).$ Then $c_1 \notin span \{ c_2 \},$ 
%and $Im(C)=span\{c_1,c_2\}.$ But this means that 
%$$Im(C) \subseteq (Im(B)+Im(D))=Im(B) \bigoplus Im(D).$$
%Similarly, $$Im(A) \subseteq Im(B) \bigoplus Im(D).$$
%Thus,
%$$dim (Im(A)+Im(B)+Im(C)+Im(D))=dim(Im(B)\bigoplus Im(D))= 4.$$

%Again, (see corollary \ref{cor:dim}), arrange  the vertices
%of the graph in a sequence  $$A, B, C, D, A_{i_5}, A_{i_6},\dots , 
%A_{i_{n-1}},$$ such that each term 
%$A_{i_j}, \,\,\, 5 \leq j \leq n-1,$
 %is 
%a friend of one the  terms 
%$$A,B,C,D, A_{i_5}, A_{i_6},\dots , A_{i_{j-1}}.$$ 
%Then, as in corollary \ref{cor:dim},
%$$\dim(Im(A)+Im(B)+Im(C)+Im(D)+Im(A_{i_5})+ \dots +Im(A_{i_{n-1}}))
%\leq 4+1\cdot(n-5) = n-1,$$ contradiction.
%So, $B$ and $D$ are friends. For the same reason $A$ and $C$ are 
%also friends.
%Thus, the 4-gon has both diagonals.

%\vskip .5cm

\begin {T} \label {T:discon} Let $\rho:B_n \to GL_r(\C)$ be an irreducible
 representation,  where $r \geq n$ and $n \geq 6.$
Let $rank (A_1)=2.$

Then  $Im(A_i) \cap Im(A_{i+1}) \neq \{0\}$ for $1 \leq i \leq n-2;$ 
that is the friendship graph of $\rho$ contains the chain graph.

\end{T}
\pf Suppose not. Then by Theorem \ref{T:dva} $(b),$
$Im(A_i) \cap Im(A_j) \neq 0$ whenever $A_i$ and $A_j$ are
not neighbors. Consider $$U=Im(A_1) +Im(A_2)+Im(A_3).$$
Since $Im(A_1) \cap Im(A_3) \neq 0,$ $dim U \leq 5.$

For $i =4, \dots ,n-1,$ let $a_i,\,\,\,b_i$ be, respectively, 
nonzero elements of $Im(A_1) \cap Im(A_i)$ and 
$Im(A_2) \cap Im(A_i).$ Since $Im(A_1) \cap Im(A_2)=0,$ 
$a_i$ and $b_i$ are linearly independent, so they are
a basis for $Im(A_i),$ and $Im(A_i) \subseteq Im(A_1) +Im(A_2).$ 
Thus $$U =Im(A_1) +Im(A_2)+\dots +Im(A_{n-1}),$$ which is invariant
under $\rho(B_n).$ Thus $r \leq 5,$ by the irreducibility of $\rho,$
a contradiction with $r \geq n \geq 6.$

%By corollary \ref{cor:connect} the friendship graph is  
%connected. Then there exists $i$ such that $A_1$ and $A_i$ are 
%friends.
%Suppose that no neighbors are friends. Then we have that by the 
%lemma \ref{lem:fr}
%about friends  any $A_k$ for $k=3,4,5, \dots, n-1$ is a friend
%of $A_1:$
%\vskip 1cm
%\psfig{figure=disc6.eps}
%\vskip 1cm

%Then, $A_2$ and $A_k$ are friends  for $k=4,5, \dots ,n-1,$ in 
%particular,
%$A_2$ and $A_4$ are friends and $A_2$ and $A_5$ are friends (here we %need
%$n \geq 6$):
%\vskip 1cm
%\psfig{figure=disc6a.eps}
%\vskip 1cm

%So, we have a 4-gon $A_1 - A_4 -A_2 - A_5$ without the diagonal
%$A_4 -A_5,$ which contradicts  lemma \ref{lem:4gon}.

\begin{Rem}\label{Rem:5} For $n=5$ and $\rho$ satisfying the hypothesis of theorem \ref{T:discon} there are two possible
friendship graphs: 1) all neighbors are friends and 2) an exceptional case:
 \vskip 1cm
\psfig{figure=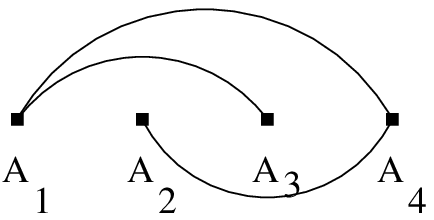}
\vskip 1cm

By \cite{thesis}, Theorem 7.1, part 2, every irreducible representation 
with the above friendship graph is equivalent to the restriction to $B_5$
of the Jones' representation (see \cite{form}, p. 296).
\end{Rem}

%%%%%%%%%%%%%%%%%%%%%%%%%%%%%%%%%%%%%%%%%%%%%%%%%%
%%%%%%%%%%%%%%%%%%%%%%%%%%%%%%%%%%%%%%%%%%%%%%%%%%%

\begin{lem}\label{lem:chain} Let $\rho:B_n \to GL_r(\C)$ be an irreducible
 representation,  where $r \geq n,$  $n \geq 5,$
and $rank (A_1)=2.$ Suppose that the associated friendship graph contains the chain.

Then $r=n$ and the associated friendship graph is the chain
(that is, the only edges are between neighbors).
\end{lem}
\pf  By corollary \ref{cor:dim}, $r=n.$ Consider the full friendship graph of  $\rho.$ Then
$$Im(A_i) \cap Im(A_{i+1}) \neq \{0\}$$ for any $i$
where indices are taken modulo $n.$ If $Im(A_i) \cap Im(A_{i+1})$
is two-dimensional, then $Im(A_1) =Im(A_2) = \dots ,$ and  $Im(A_1)$ is a two-dimensional invariant subspace, contradicting the irreducibility of $\rho.$ 
Hence $Im(A_i) \cap Im(A_{i+1})$ are one-dimensional.

For any $x \in Im(A_i),$ $x=A_iy,$ $x \neq 0,$ we have that $$ T x = T A_i y = T A_i T ^{-1}
(T y)=A_{i+1} (T y) \in Im(A_{i+1})$$  for $T= \rho(\tau).$ Moreover, $T x \neq 0$ because
$T$ is invertible.

Choose $x_1 \neq 0$ to be a basis vector for $Im(A_1) \cap Im(A_2).$
Define $x_{i+1}= T ^{i} x_1$ for $1\leq i \leq n-1.$
Then $x_{i}$ is a basis vector for $Im(A_i) \cap Im(A_{i+1}).$

If for some $i,$ $x_i$ is proportional to $x_{i+1}$ then, because
a full friendship graph is a ${\Z}_n$-graph, all the $x_j$ are proportional
to $x_1.$ Then, because we have 5 or more vertices in the
full friendship graph, for any $A_i$ there exists $j$ such that
both $A_j$ and $A_{j+1}$ are not neighbors of $A_i.$ Then
$$A_i A_j=A_j A_i$$ and $$A_i A_{j+1}=A_{j+1}A_i.$$ So, if $x \in Im(A_j) \cap Im(A_{j+1})$ then $A_i x \in  Im(A_j) \cap Im(A_{j+1}).$
But this means that $ span \{x_1\}$ is an invariant subspace and
the representation is not irreducible. 

So, if the representation is irreducible, then for any $i,$
$x_i \notin span \{x_{i+1}\}.$ From this follows that  for any $i$
$$ Im(A_i)= span \{x_{i-1}, x_i\}$$ and the $n$ vectors $x_0, x_1, \dots, x_{n-1}$
form a basis of $V.$
Then for any two non-neighbors $A_i$ and $A_j$ $$Im(A_i) \cap Im(A_{j})=\{0\}.$$

%\vskip .5cm

Now, we have the following
\begin{T}\label{T:4} Let $\rho: B_n \to GL_r ({\C})$ be irreducible, where $r \geq n.$ Suppose that for any generator $\sig_i,$
$\rho (\sig_i)=1+A_i,$ where $rank (A_i)=2.$

1) If $n \geq 6,$ then $r=n$ and $\rho$ has a friendship graph which is a chain.

2) If $n=5,$ then  $r=5$ and either $\rho$ has a friendship graph which is a chain or
$\rho$ has the exceptional friendship graph {\rm (see Remark \ref{Rem:5}).}

3) If $n=4,$ then either $r=4$ and $\rho$ has a  friendship graph which is a chain; or
$\rho$ has one of the following  exceptional friendship graphs:

\vskip 1cm
\psfig{figure=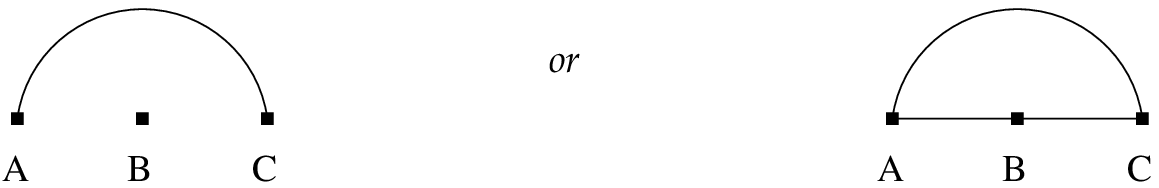}

\vskip 1cm

\end{T} 
\pf 1) If $n \geq 6,$ then by theorem \ref{T:discon} the associated friendship
graph contains a chain, and, by lemma \ref{lem:chain} has no other edges 
and $r=n.$

2) If $n=5,$ then by corollaries \ref{cor:connect} and \ref{cor:dim}
the friendship graph of $\rho$ is connected and $r=n.$ If it  contains a chain
graph, then, by lemma \ref{lem:chain}, it has no other edges. If it does not
contain a chain graph, we obtain the exceptional case.

3) If $n =4,$ then by theorem \ref{T:38} the friendship graph is not
totally disconnected. Hence, we have only three possible ${\Z}_4-$graphs 
on 4 vertices.

\begin{Rem} It is proven in \cite{thesis}, Chapter 6, that
any representation of $B_4$ with either of the exeptional friendship
graphs in 3) of the above theorem is reducible.

\end{Rem}
%%%%%%%%%%%%%%%%%%%%%%%%%%%%%%%%%%%%%%%%%%%%%%%%%%%%%%%%%%%%
\section{Representations whose friendship graph is a chain}

\begin{defi} {\bf The standard representation} is the representation
$$\tau_n:B_n \to GL_n ({\Z}[t ^{\pm 1}]$$ defined by

\vskip 1cm
$$\tau_n (\sig_i)=\left( \begin{array}{ccccc}
I_{i-1}&&&\\
&0&t&\\
&1&0&\\
&&&I_{n-1-i}
\end{array} \right),$$ 
\vskip .5cm
\noindent
for $i=1,2,\dots, n-1,$ where $I_{k}$ is the $k\times k$ identity matrix.

\end{defi}

\begin{T} \label{T:cha} Let $\rho : B_{n} \to GL_{n}({\C})$ be an irreducible representation, where $n \geq 4.$
Suppose that $\rho(\sig_1)=1+A_1,$ where $rank (A_1)=2,$ and the associated friendship graph of $\rho$ is a chain.

Then $\rho$ is equivalent to a specialization $\tau_n(u)$ of the standard representation for some $u \in {\C} ^*.$

\vskip 1cm
\psfig{figure=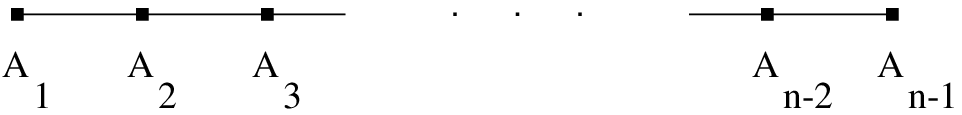}
\vskip 1cm

\end{T}

Before proving the theorem, we will need the following technical lemma:

\begin{lem} \label{lem:basis} Let $A$ be a friend and a neighbor of $B,$
$B$ be a friend and a neighbor of $C$ and suppose that $A$  is not a friend of $C$ :
\vskip 1cm
\psfig{figure=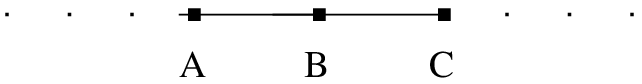}
\vskip 1cm
Let $a \neq 0$ be such that $span\{a\}= Im(A) \cap Im(B),$
and let $b=(1+B)a.$ %and $c=(1+C)b.$ 
Then: 

1) $span\{b\} = Im (C)\cap Im(B).$

%2) $Bc=0.$

2) $%(1+C)c \in span\{b\},\,\,\,
 (1+B)b \in span\{a\}$ and $(1+B)b \neq 0.$
 
3) The vectors $a$ and $b$ are linearly independent.

\end{lem}
\pf First of all, notice that the vector $b$ is non-zero, because
$1+B$ is invertible and $a \neq 0.$

1) $b=(1+B)a \in Im(B),$ because $a \in Im(B).$

$A$ and $C$ are not friends, that is $CA=0,$ so $Ca=0.$ Let $a =Ba_1.$ Then
$$(1+B)a=(1+B+BC)a=(1+B+BC)Ba_1=(B+B^2+BCB)a_1=$$
$$=(C+C^2+CBC)a_1 \in Im(C);$$
that is, $b \in Im (C)\cap Im(B),$ and because $ Im (C)\cap Im(B) $ is
one-dimensional and $b \neq 0,$
 $$span\{b\} = Im (C)\cap Im(B).$$ 

%2) $$Bc=B(1+C)b=B(1+C)(1+B)a=B(1+B+C+CB)a=$$
%$$=(B+B^2+BCB)a+BCa=(C+C^2+CBC)a+BCa=$$
%$$=(1+C+CB+B)Ca=0$$
%because $Ca=0.$

2) Clearly, $(1+B)b \in Im(B).$

Note, that $Ab=0,$ as $b \in Im(C)$ by the above, and $AC=0.$ Let $b=Ba^{'}.$
Then $$(1+B)b=(1+B+BA)b=(1+B+BA)Ba^{'}=(A+A^2+ABA)a^{'} \in Im(A).$$

%To prove this statement we need to show that $(1+C)c \in Im(C) \cap Im(B).$ %Indeed,
%$(1+C)c \in Im(C)$ because $c \in Im(C);$ and because  $Bc=0$, 
%$$(1+C)c=(1+C+CB)c=(1+C+CB)Cc_1=(C+C^2+CBC)c_1=$$
%$$(B+B^2+BCB)c_1 \in Im(B),$$ where $c=Cc_1,$ so $(1+C)c \in span\{b\}.$

 %The fact that  $(1+B)b \in span\{a\}$ can be proven in the same way
%noticing that $Ab=0$ (cf. corollary \ref{cor:ABC}, 2).)

3) $a \in Im(A),$  $b \in Im(C)$ by part 1), and $ Im (A)\cap Im(C)=\{0\}$
by the hypothesis of the lemma.

\vskip .5cm

%\begin{cor} \label {cor:ABC} Retaining the hypothesis of lemma \ref{lem:basis},

%1) $a$ and $b$ are linearly independent, and $Im(B)=span\{a,b\}.$

%2) $Ab=0.$

%3) $b$ and $c$ are linearly independent, and $Im(C)= span\{b,c\}.$

%4) $ Im (A)\cap Im(B)=span\{(1+B)b\},$

 %$ Im (B)\cap Im(C)=span\{(1+C)c\}.$

%\end{cor}
%\pf 1) $b \notin span\{a\},$ because $ Im (C)\cap Im(A)=\{0\}$ and $a \in %Im(A),$  $b \in Im(C).$

%2) $ AC=0,$  $b \in Im(C),$ hence  $Ab=0.$

%3) If $b$ and $c$ are linearly dependent then $c= \lam b$ with some
%$\lam \neq 0.$ Then $Bc=\lam Bb.$ Together with $Bc=0$ this means that
%$Bb=0.$ But by the lemma above, part 3)  $(1+B)b \in span\{a\},$
%and if $Bb=0,$ then $b \in span\{a\},$ which contradicts 
%lemma  \ref{lem:basis}, part 1).

%4) These follow from the above and the fact that $(1+B)b$ and $(1+C)c$ are
%non-zero.

%\begin{Rem} All the statements above are true even if $C$ 
%is the last  vertex in the chain.
%\end{Rem}

\vskip .5cm
{\bf Proof of Theorem \ref{T:cha}} We include the redundant generator $\sig_0,$
and indices are modulo $n.$  Consider  $Im(A_i) \cap Im(A_{i+1}),$
which is $0,$ $1,$ or $2-$dimensional. It is nonzero, because of the hypothesis
that the friendship graph is a chain. It is not 2-dimensional, for then
$$Im(A_0)=Im(A_1) = \dots =Im(A_{n-1})$$ would be a $2-$dimensional invariant subspace, contradicting the irreducibility of $\rho.$ Hence,  $Im(A_i) \cap Im(A_{i+1})$ is one-dimensional.

%%%%%%%%%%%%%%%%%%%%%%%%%%%%%%%%%%

Let $a_0$ be a basis vector for $Im(A_0) \cap Im(A_{1}).$ Let 
$$a_1=(1+A_1)a_0,\,\,\,a_2=(1+A_2)a_1,\,\,\,\dots ,\,\,\,a_{n-1}=(1+A_{n-1})a_{n-2}.$$ By induction and lemma \ref{lem:basis},
part 1), $a_i$ is a basis vector for  $Im(A_i) \cap Im(A_{i+1}),$
for  $0 \leq i \leq n-1.$ By lemma \ref{lem:basis}, part 3), $a_i$ and $a_{i+1}$
are linearly independent. Thus $\{a_i,a_{i+1} \}$ is a basis for $Im(A_i).$

Since $$span \{ a_0, \dots a_{n-1} \}=Im(A_1)+ \dots +Im(A_{n-1})$$ 
is invariant under $B_n$ and $\rho$ is an $n-$dimensional
irreducible representation, $\{ a_0, \dots a_{n-1} \}$ is a basis
for ${\C}^n.$

We now wish to determine the action of 
$\rho(\sig_1), \,\,\,\rho(\sig_2), \dots , \rho(\sig_{n-1})$ on this basis.

Consider $a_i \in Im(A_i) \cap Im(A_{i+1}).$ If $j \neq i,\,\,\,i+1,$
then $A_j$ is not a neighbor of one of $A_i,\,\,\,A_{i+1}$ (since $n \geq 4$),
say $A_k,$ and then $A_kA_j=A_jA_k=0,$ so $A_ja_i=0,$ and 
$$\rho(\sig_j)a_i=(1+A_j)a_i=a_i.$$

By our construction $$\rho(\sig_{i+1})a_{i}=(1+A_{i+1})a_{i}=a_{i+1}$$ for 
$0 \leq i \leq n-2.$

By lemma \ref{lem:basis}, part 2), $$\rho(\sig_i)a_i=(1+A_i)a_i=u_ia_{i-1},$$
for $1 \leq i \leq n-1,$ where $u_i \in {\C}^*.$ 

By the above calculations the matrices of $\rho(\sig_1), \dots ,\rho(\sig_{n-1})$ with respect to the basis $a_0,\,\,\,a_1, \dots ,a_{n-1}$
are

%%%%%%%%%%%%%%%%%%%%%%%%%%%%%%%%%%%%%%%%%%%%%%%%%%%%%%%%%

%Consider a vector $a'_2 \neq 0$ such that $span\{a'_2\}=Im(A_2) \cap Im(A_3).$
%Choose  vectors $a'_1=(1+A_2)a'_2$ and $a_0=(1+A_1)a'_1.$ By corollary
%\ref{cor:ABC}, part 3), $a_0$ and $a'_1$ are linearly independent,
%$Im(A_1)=span\{a_0,a'_1\},$ and by lemma \ref{lem:basis}, part 1),
%vector  $a'_1$ spans $Im(A_1) \cap Im(A_2).$ 
%Now consider vectors $$a_1=(1+A_1)a_0,\,\,\,a_2=(1+A_2)a_1,\,\,\, \dots
%a_{i+1}=(1+A_{i+1})a_i,\,\,\,a_{n-1}=(1+A_{n-1})a_{n-2}.$$ 
%Then by lemma \ref{lem:basis}, part 3), $$span\{a_1\}=span\{a'_1\},$$
% $$span\{a_2\}=span\{a'_2\},$$ and  by lemma \ref{lem:basis}, part 1),
%for every $i,\,\,1 \leq n-2$
%$$span\{a_i\}=Im(A_{i}) \cap Im(A_{i+1})$$ and 
%$$Im(A_{i})=span \{a_{i-1},a_{i}\}.$$
%Also, by corollary \ref{cor:ABC}, part 3) we have that 
%$$Im(A_{n-1})=span \{a_{n-2},a_{n-1}\}.$$

%Now,  $$V=Im(A_1)+Im(A_2)+ \dots +Im(A_{n-1})=span\{a_0,a_1, \dots ,a_{n-1}\}$$
%and $dim V=n,$ so vectors $a_0,a_1, \dots ,a_{n-1}$ form a basis of $V.$
%%%%%%%%%%%%%%%%%%%%%%%%%%%%%%%%%%%%%%%%%%%%%5

\vskip 1cm
$$\rho (\sig_i)=\left( \begin{array}{ccccc}
I_{i-1}&&&\\
&0&u_i&\\
&1&0&\\
&&&I_{n-1-i}
\end{array} \right),$$ 
\vskip .5cm
\noindent
for $i=1,2,\dots, n-1,$ where $I_{k}$ is the $k\times k$ identity matrix,
and $u_1, \dots ,u_{n-1} \in \C^{\ast}.$ Since $\sig_1 ,\dots ,\sig_{n-1}$
are conjugate in $B_n,$ the $u_i$ are all equal, and we have 
the standard representation.
\vskip .5cm

Now let us consider when the standard representation is irreducible.
\begin{lem} \label{lem:54} If $u=1$ then  $\tau_n(u)$ is reducible.

\end{lem}
\pf If $u=1$ then the vector $v=(1,1,1, \dots , 1)^{T}$ is a fixed vector. 

\vskip .5cm

\begin{lem} \label{lem:55} If $u \neq 1$ then  $\tau_n(u)$  is irreducible.

\end{lem}
\pf We need to prove that starting from any non-zero vector
$x=\sum a_i e_i,$ we can generate the whole space. Obviously,
it is enough to show that we can generate one of the standard
basis vectors $e_i$. To do this, take $i$ such that $a_i \neq 0.$
Consider the operator
$$H=A+A^2+ABA=B+B^2+BAB,$$
where $A=\rho(\sig_{i-1})$ and $B=\rho(\sig_i)$. By a direct
calculation $Hx=(u-1)a_ie_i$. Because $u\neq 1$ the vector $Hx$
is a non-zero multiple of $e_i$.

\vskip .5cm

Now, we have the main result of this paper:

\begin{T} [The Main Theorem]  \label {T:main}
Let $\rho : B_n \to GL_r ({\C})$ be an irreducible representation
of $B_n$ for $n \geq 6.$ Let $r \geq n,$ and let
$\rho (\sig_1)=1+A_1$  with
$rank (A_1)=2.$ 

Then $r=n$ and $\rho$ is equivalent to the following
representation :
$$\tau:B_n \to  GL_n ({\C}),$$

\vskip 1cm
$$\rho (\sig_i)=\left( \begin{array}{ccccc}
I_{i-1}&&&\\
&0&u&\\
&1&0&\\
&&&I_{n-1-i}
\end{array} \right),$$ 
\vskip .5cm
\noindent
for $i=1,2,\dots, n-1,$ where $I_{k}$ is the $k\times k$ identity matrix,
and $u \in {\C}^*,$ $u \neq 1.$ These representations are non-equivalent
for different values of $u.$
\end{T}
\pf By Theorem \ref{T:4} the friendship graph of $\rho$ is a chain. Then,
by theorem \ref{T:cha}, $\rho$ is equivalent to a standard representation
$\tau (u)$ for some $u \in {\C}^*.$ By Lemmas \ref{lem:54} and \ref{lem:55}
$u \neq 1.$

Combining Theorem \ref{T:main} and the classification theorem of Formanek
(see \cite{form}, Theorem 23), we get the following

\begin{cor} Let $\rho : B_n \to GL_r ({\C})$ be an irreducible representation
of $B_n$ for $n \geq 7.$ Let $corank(\rho)=2.$  

Then $\rho$ is equivalent to a specialization of the standard representation
$\tau_n(u),$ for some $u \neq 1, \,\,\, u \in {\C}^*.$ 

\end{cor}

%%%%%%%%%%%%%%%%%%%%%%%%%%%%%%%%%%%%%%%%%%%%%%%%%%%%%%%%%%


\begin{thebibliography}{99}
\bibitem{birm} J.S. Birman, {\it Braids, Links, and Mapping Class Groups,}
Ann.of Math.Stud.,82, Princeton Univ.Press,Princeton,N.J.,1974.

\bibitem{chow} W.-L. Chow, {\it On the algebraical braid group,} Ann. of Math.
{\bf 49} (1948), 654-658.

\bibitem{form} E. Formanek, {\it Braid Group Representations of Low Degree,}
Proc.London Math.Soc. {\bf 73} (1996), 279-322.

%\bibitem{han} V.L. Hansen, {\it Braids and Coverings. Selected topics.}
%Cambridge Univ. Press ,1989.

\bibitem{jones} V.F.R.Jones, {\it Hecke algebra representations of braid groups
and link polynomials,} Ann. of Math.{\bf 126} (1987), 335-388.


%\bibitem{lee}  W.Lee, {\it Representations of the braid group
%$B_4,$} J.Korean Math.Soc. {\bf 34} (1997), no. 3, 673-693.

\bibitem{thesis} I. Sysoeva, {\it On the irreducible representations
of braid groups,}  Ph.D. thesis,  1999.

\bibitem{tym} Dian-Min Tong, Shan-De Yang, Zhong-Qi Ma, {\it A new class of representations of braid groups,}  Comm. Theoret. Phys. {\bf 26} (1996), no. 4, 483--486.





\end{thebibliography}
\end{document}